\documentclass[10pt]{article}
\usepackage{cite}
\usepackage{mathrsfs,textcomp}
\usepackage{amsthm}
\usepackage{amssymb,amsmath}
\newtheorem{defn}{Definition}[section]

\newtheorem{thm}{Theorem}[section]

\newtheorem{prop}[thm]{Proposition}

\newtheorem{example}[thm]{Example}

\usepackage{enumitem}
\begin{document}
	\title{\textbf{Solving system of Urysohn type integral equation via bicomplex partial b-metric space }}
	\author{Gunaseelan Mani$^{1}$, Arul Joseph Gnanaprakasam$^{2}$, Yongjin Li $^{3}$,\\ and Zhaohui Gu $^{4,}$* \\}
	\date{	
		$^{1}$Department of Mathematics,\\
		Sri Sankara Arts and Science College(Autonomous),\\
		Affiliated to Madras University,\\
		Enathur, Kanchipuram, Tamil Nadu, India 631 561.\\
		E-mail: mathsguna@yahoo.com\\
		$^{2}$Department of Mathematics,\\
		College of Engineering and Technology,\\
		Faculty of Engineering and Technology,\\
		SRM Institute of Science and Technology,\\
		SRM Nagar, Kattankulathur 603203,\\
		Kanchipuram, Chennai, Tamil Nadu, India.\\
		E-mail: aruljoseph.alex@gmail.com\\ 
		$^{3}$Department of Mathematics,\\
		Sun Yat-Sen University,\\
		Guangzhou 510275, China.\\
		E-mail: stslyj@mail.sysu.edu.cn\\
		and\\
		$^{4}$School of Mathematics and Statistics,\\
		Guangdong University of Foreign Studies,\\
		Guangzhou 510006, China.\\
		E-mail: zhgugz@163.com\\
	} 
	\maketitle 
	\begin{abstract}
		\textbf{In this paper, we introduce the notion of bicomplex partial b-metric space and prove some common fixed point theorems. Our results generalize and expand some of the literature's well-known results. We also explore some of the applications of our key results to Urysohn type integral equations.}
	\end{abstract}
	\section{Introduction}
	Serge \cite{bc1} made a pioneering attempt in the development of special algebas. He conceptualized commutative generalization of complex numbers as bicomplex numbers, tricomplex numbers, etc. as elements of an infinite set of algebra. Subsequently during the $1930$, other researchers also contributed in this area \cite{bc2}-\cite{bc8}. But unfortunately the next fifty years failed to witness any advancement in this field. Afterward Price \cite{bc6} developed the bicomplex algebra and function theory. Recently renewed interest in this subject finds some significant applications in different fields of mathematical sciences as well as other branches of science and technology. Also one can see the attempts in \cite{bc9}. An impressive body of work has been developed by a number of researchers. Among them an important work on elementary functions of bicomplex numbers has been done by Luna-Elizaarrar$\acute{a}$s, Shapiro, Struppa and Vajiac \cite{bc5}. Choi , Datta, Biswa and Islam \cite{bc3} proved some common fixed point theorems in connection with two weakly compatible mappings in bicomplex valued metric spaces. Jebril \cite{bc4} proved some common fixed point
	theorems under rational contractions for a pair of mappings in bicomplex valued metric spaces. In 2021, Beg, Kumar Datta and Pal \cite{bc10} proved fixed point theorems on bicomplex valued metric spaces. In 2021, Datta, Pal, Sarkar and Manna \cite{bc12} proved common fixed point theorem in bicomplex valued b-metric
	space.
	Gunaseelan \cite{1aa} presented the concept of complex partial b-metric space in 2019, as well as proving the fixed point theorem under the contractive condition. Many researchers have studied some intriguing concepts and applications and  has  shown significant results [\citen{bc13}-\citen{bc14}]. In this paper, we prove some common fixed point theorems on bicomplex partial b-metric space.
	\section{Preliminaries}
	Throughtout this paper, we denote the set of real, complex and bicomplex numbers respectively as $\mathbb{C}_{0}$, $\mathbb{C}_{1}$ and $ \mathbb{C}_{2}$ Segre \cite{bc1} defined the bicomplex number as:
	\begin{align*}
	\sigma=\mathfrak{a}_{1}+\mathfrak{a}_{2}i_{1}+\mathfrak{a}_{3}i_{2}+\mathfrak{a}_{4}i_{1}i_{2},
	\end{align*}
	where $\mathfrak{a}_{1},\mathfrak{a}_{2},\mathfrak{a}_{3},\mathfrak{a}_{4}\in \mathbb{C}_{0}$, and independent units $i_{1},i_{2}$ are such that $i_{1}^{2}=i_{2}^{2}=-1$ and $i_{1}i_{2}=i_{2}i_{1}$,  we denote the set of bicomplex numbers $ \mathbb{C}_{2}$ is defined as:
	\begin{align*}
	\mathbb{C}_{2}=\{\sigma:\sigma=\mathfrak{a}_{1}+\mathfrak{a}_{2}i_{1}+\mathfrak{a}_{3}i_{2}+\mathfrak{a}_{4}i_{1}i_{2},\mathfrak{a}_{1},\mathfrak{a}_{2},\mathfrak{a}_{3},\mathfrak{a}_{4}\in \mathbb{C}_{0}\},
	\end{align*}
	i.e.,
	\begin{align*}
	\mathbb{C}_{2}=\{\sigma:\sigma=\varphi_{1}+i_{2}\varphi_{2},\varphi_{1},\varphi_{2}\in \mathbb{C}_{1}\},
	\end{align*}
	where $\varphi_{1}=\mathfrak{a}_{1}+\mathfrak{a}_{2}i_{1}\in \mathbb{C}_{1}$ and $\varphi_{2}=\mathfrak{a}_{3}+\mathfrak{a}_{4}i_{1}\in \mathbb{C}_{1}$.
	If $\sigma=\varphi_{1}+i_{2}\varphi_{2}$ and $\vartheta=\mathfrak{w}_{1}+i_{2}\mathfrak{w}_{2}$ be any two bicomplex numbers then the sum is $\sigma\pm\vartheta=(\varphi_{1}+i_{2}\varphi_{2})\pm(\mathfrak{w}_{1}+i_{2}\mathfrak{w}_{2})=\varphi_{1}\pm\mathfrak{w}_{1}+i_{2}(\varphi_{2}\pm\mathfrak{w}_{2})$ and the product is $\sigma.\vartheta=(\varphi_{1}+i_{2}\varphi_{2})(\mathfrak{w}_{1}+i_{2}\mathfrak{w}_{2})=(\varphi_{1}\mathfrak{w}_{1}-\varphi_{2}\mathfrak{w}_{2})+i_{2}(\varphi_{1}\mathfrak{w}_{2}+\varphi_{2}\mathfrak{w}_{1})$. \\
	There are four idempotent elements in $ \mathbb{C}_{2}$, they are $0,1,\mathfrak{e}_{1}=\frac{1+i_{1}i_{2}}{2},\mathfrak{e}_{2}=\frac{1-i_{1}i_{2}}{2}$ out of which $\mathfrak{e}_{1}$ and $\mathfrak{e}_{2}$ are nontrivial such that $\mathfrak{e}_{1}+\mathfrak{e}_{2}=1$ and $\mathfrak{e}_{1}\mathfrak{e}_{2}=0$. Every bicomplex number $\varphi_{1}+i_{2}\varphi_{2}$ can be uniquely be expressed as the combination of $\mathfrak{e}_{1}$ and $\mathfrak{e}_{2}$, nemely
	\begin{align*}
	\sigma=\varphi_{1}+i_{2}\varphi_{2}=(\varphi_{1}-i_{1}\varphi_{2})\mathfrak{e}_{1}+(\varphi_{1}+i_{1}\varphi_{2})\mathfrak{e}_{2}.
	\end{align*}
	This representation of $\sigma$ is known as the idempotent representation of bicomplex number and the complex coefficients $\sigma_{1}=(\varphi_{1}-i_{1}\varphi_{2})$ and $\sigma_{2}=(\varphi_{1}+i_{1}\varphi_{2})$ are known as idempotent components of the bicomplex number $\sigma$.\\
	An element $\sigma=\varphi_{1}+i_{2}\varphi_{2}\in  \mathbb{C}_{2}$ is said to be invertible if there exists another element $\vartheta$ in $ \mathbb{C}_{2}$ such that $\sigma\vartheta=1$ and $\vartheta$ is said to be inverse(multiplicative) of $\sigma$. Consequently $\sigma$ is said to be the inverse(multiplicative) of $\vartheta$. An element which has an inverse in $ \mathbb{C}_{2}$ is said to be the nonsingular element of $ \mathbb{C}_{2}$ and an element which does not have an inverse in $ \mathbb{C}_{2}$ is said to be the singular element of $ \mathbb{C}_{2}$.\\
	An element $\sigma=\varphi_{1}+i_{2}\varphi_{2}\in  \mathbb{C}_{2}$ is nonsingular if and only if $|\varphi_{1}^{2}+\varphi_{2}^{2}|\neq 0$ and singular if and only if $|\varphi_{1}^{2}+\varphi_{2}^{2}|=0$
	The inverse of $\sigma$ is defined as
	\begin{align*}
	\sigma^{-1}=\vartheta=\frac{\varphi-i_{2}\varphi_{2}}{\varphi_{1}^{2}+\varphi_{2}^{2}}
	\end{align*} 
	Zero is the only element in $\mathbb{C}_{0}$ which does not have multiplicative inverse and in $\mathbb{C}_{1}$, $0=0+i0$ is the only element which does not have multiplicative inverse. We denote the set of singular elements of $\mathbb{C}_{0}$ and $\mathbb{C}_{1}$ by $\mathfrak{O}_{0}$ and $\mathfrak{O}_{1}$ respectively. But there are more than one element in $ \mathbb{C}_{2}$ which do not have multiplicative inverse, we denote this set by $\mathfrak{O}_{2}$ and clearly $\mathfrak{O}_{0}=\mathfrak{O}_{1}\subset \mathfrak{O}_{2}$.\\
	A bicomplex number $\sigma=\mathfrak{a}_{1}+\mathfrak{a}_{2}i_{1}+\mathfrak{a}_{3}i_{2}+\mathfrak{a}_{4}i_{1}i_{2}\in  \mathbb{C}_{2}$ is said to be degenerated if the matrix 
	\begin{gather*}
	\begin{pmatrix}
	\mathfrak{a}_{1} & \mathfrak{a}_{2} \\
	\mathfrak{a}_{3} & \mathfrak{a}_{4}
	\end{pmatrix}
	\end{gather*}
	is degenerated. In that case $\sigma^{-1}$ exists and it is also degenerated.\\
	The norm $||.||$ of $ \mathbb{C}_{2}$ is a positive real valued function and $||.||: \mathbb{C}_{2}\to \mathbb{C}_{0}^{+}$ is defined by 
	\begin{align*}
	||\sigma||&=||\varphi_{1}+i_{2}z_{2}||=\{|\varphi|^{2}+|\varphi|^{2}\}^{\frac{1}{2}}\\
	&=\bigg[\frac{|(\varphi_{1}-i_{1}\varphi_{2})|^{2}+|(\varphi_{1}+i_{1}\varphi_{2})|^{2}}{2}\bigg]^{\frac{1}{2}}\\
	&=(\mathfrak{a}_{1}^{2}+\mathfrak{a}_{2}^{2}+\mathfrak{a}_{1}^{2}+\mathfrak{a}_{3}^{2}+\mathfrak{a}_{4}^{2})^{\frac{1}{2}},
	\end{align*}  
	where $\sigma=\mathfrak{a}_{1}+\mathfrak{a}_{2}i_{1}+\mathfrak{a}_{3}i_{2}+\mathfrak{a}_{4}i_{1}i_{2}=\varphi_{1}+i_{2}\varphi_{2}\in  \mathbb{C}_{2}$.\\
	The linear space $ \mathbb{C}_{2}$ with respect to defined norm is a norm linear space, also $ \mathbb{C}_{2}$ is complete, therefore $ \mathbb{C}_{2}$ is the  Banach space. If $\sigma,\vartheta\in  \mathbb{C}_{2}$ then $||\sigma\vartheta||\leq \sqrt{2}||\sigma|| ||\vartheta||$ holds instead of $||\sigma\vartheta||\leq||\sigma|| ||\vartheta||$,
	therfore $ \mathbb{C}_{2}$ is not the Banach algebra.
	The partial order relation $\preceq_{i_{2}}$ on $ \mathbb{C}_{2}$ is defined as:
	Let $ \mathbb{C}_{2}$ be the set of bicomplex numbers and $\sigma=\varphi_{1}+i_{2}\varphi_{2}$, $\vartheta=\mathfrak{w}_{1}+i_{2}\mathfrak{w}_{2}\in  \mathbb{C}_{2}$ then $\sigma\preceq_{i_{2}} \vartheta$ if and only if $\varphi_{1}\preceq \mathfrak{w}_{1}$ and $\varphi_{2}\preceq\mathfrak{w}_{2}$, i.e., $\sigma\preceq_{i_{2}} \vartheta$ if one of the following conditions is satisfied:
	\begin{enumerate}[label=(\Alph*)]
		\item $\varphi_{1}=\mathfrak{w}_{1}$, $\varphi_{2}=\mathfrak{w}_{2}$,
		\item $\varphi_{1}\prec\mathfrak{w}_{1}$, $\varphi_{2}=\mathfrak{w}_{2}$,\label{be1}
		\item $\varphi_{1}=\mathfrak{w}_{1}$, $\varphi_{2}\prec\mathfrak{w}_{2}$, $\text{and}$\label{be2}
		\item  $\varphi_{1}\prec\mathfrak{w}_{1}$, $\varphi_{2}\prec\mathfrak{w}_{2}$,\label{be3}
	\end{enumerate}
	In particular we can write $\sigma\lnsim_{i_{2}} \vartheta$ if $\sigma\preceq_{i_{2}} \vartheta$ and $\sigma\neq \vartheta$ i.e., one of \ref{be1},\ref{be2} and \ref{be3} is satisfied and we will write $\sigma\prec_{i_{2}}\vartheta$ if only \ref{be3} is satisfied.\\
	For any two bicomplex numbers $\sigma,\vartheta\in  \mathbb{C}_{2}$ we can verify the followings:
	\begin{enumerate}[label=(S\arabic*)]
		\item $\sigma\preceq_{i_{2}} \vartheta \Rightarrow ||\sigma||\leq ||\vartheta||$,
		\item $||\sigma+\vartheta||\leq ||\sigma||+||\vartheta||$,
		\item  $||\mathfrak{a}\sigma||=\mathfrak{a}||\sigma||$, where $\mathfrak{a}$ is a non negative real number,
		\item  $||\sigma\vartheta||\leq \sqrt{2}||\sigma||||\vartheta||$ and the equality holds only when at least one of $\sigma$ and $\vartheta$ is degenerated,
		\item  $||\sigma^{-1}||=||\sigma||^{-1}$ if $\sigma$ is a degenerated bicomplex number with $0\prec \sigma$,
		\item $||\frac{\sigma}{\vartheta}||=\frac{||\sigma||}{||\vartheta||}$, if $\vartheta$ is a degenerated bicomplex number. 
	\end{enumerate}
	Now, let us recall some basic concepts and notations, which will be used in the sequel.

	\begin{defn}\label{dt1}
		A bicomplex partial b-metric on a non-void set $\OE$ is a function $\delta_{cb}:\OE \times \OE \to \mathbb{C}^{+}_{2}$ such that for all $\sigma,\vartheta,\varphi\in \OE$:\\
		(i) $0\preceq_{i_{2}} \delta_{cb}(\sigma, \sigma) \preceq_{i_{2}} \delta_{cb}(\sigma,\vartheta)(\text{small self-distances})$\\
		(ii) $\delta_{cb}(\sigma,\vartheta)=\delta_{cb}(\vartheta, \sigma) (symmetry)$\\
		(iii) $\delta_{cb}(\sigma, \sigma)=\delta_{cb}(\sigma,\vartheta)=\delta_{cb}(\vartheta,\vartheta) \Leftrightarrow \sigma=\vartheta (equality)$\\
		(iv) $\exists$ a real number $s\geq 1$ and $s$ is an independent of $\sigma,\vartheta,\varphi$ such that $\delta_{cb}(\sigma,\vartheta)\preceq_{i_{2}} s[\delta_{cb}(\sigma,\varphi)+\delta_{cb}(\varphi,\vartheta)]-\delta_{cb}(\varphi,\varphi) (triangularity)$.\\
		A bicomplex partial b-metric space is a pair $(\OE,\delta_{cb})$ such that $\OE$ is a non-void set and $\delta_{cb}$ is the bicomplex partial b-metric on $\OE$. The number $s$ is called the coefficient of $(\OE,\delta_{cb})$. 
	\end{defn}
	\begin{example}
		Let $\OE=\mathbb{R}^{+}$, $q>1$ a constant and $\delta_{cb}:\OE\times \OE \to \mathbb{C}^{+}_{2}$ be defined by 
		$\delta_{cb}(\sigma,\vartheta)=[max\{\sigma,\vartheta\}]^q+\lvert \sigma-\vartheta \rvert^q+i\{[max\{\sigma,\vartheta\}]^q+\lvert \sigma-\vartheta \rvert^q\}$ $\forall$ $\sigma,\vartheta \in \OE$. 
		Then $(\OE, \delta_{cb})$ is a bicomplex partial b-metric space with coefficient $s=2^q>1$, but it is neither a bicomplex valued b-metric nor  a bicomplex partial metric. Indeed, for any $\sigma>0$ we have $\delta_{cb}(\sigma,\sigma)=\sigma^p(1+i)\neq 0$. Therefore, $\delta_{cb}$ is not a bicomplex valued b-metric on $\OE$. Also, for $\sigma=6, \vartheta=2, \varphi=5$
		\begin{align*}
		\delta_{cb}(\sigma,\vartheta)&=(6^q+5^q)(1+i),
		\end{align*}
		\begin{align*}
		\delta_{cb}(\sigma,\vartheta)+\delta_{cb}(\vartheta,\vartheta)-\delta_{cb}(\vartheta,\vartheta)&=(6^q+1^q)(1+i)+(5^q+3^q)(1+i)-5^q(1+i)\\
		&=(6^q+1+3^q)(1+i).
		\end{align*}
		So,  $\delta_{cb}(\sigma,\vartheta)\succ_{i_{2}}\delta_{cb}(\sigma,\vartheta)+\delta_{cb}(\vartheta,\vartheta)-\delta_{cb}(\vartheta,\vartheta)$ $\forall$ $q>1$. 
		Therefore $\delta_{cb}$ is not a bicomplex partial metric on $\OE$. 
	\end{example}
	\begin{prop}
		Let $\OE$ be a non-empty set such that $\mathfrak{p}$ is a bicomplex partial metric and $\mathfrak{d}$ is a bicomplex valued b-metric with coefficient $s>1$ on $\OE$. Then the function $\delta_{cb}: \OE\times\OE \to \mathbb{C}^{+}_{2}$ defined by $\delta_{cb}(\sigma,\vartheta)=\mathfrak{p}(\sigma,\vartheta)+\mathfrak{d}(\sigma,\vartheta)$ $\forall$ $\sigma,\vartheta \in \OE$ is a bicomplex partial b-metric on $\OE$, that is, $(\OE, \delta_{cb})$ is a bicomplex partial b-metric space. 
	\end{prop}
	\begin{proof}
		Let $(\OE, \mathfrak{p})$ be a bicomplex partial metric space and $(\OE, \mathfrak{d})$ be a bicomplex b-metric space with coefficient $s>1$. Then $(i)$, $(ii)$, and $(iii)$ of Definition \ref{dt1} are obvious for the function $\delta_{cb}$. Let $\sigma, \vartheta, \varphi\in \OE$ be arbitrary.\\
		Now,
		\begin{align*}
		\delta_{cb}(\sigma,\vartheta)&=\mathfrak{p}(\sigma,\vartheta)+\mathfrak{d}(\sigma,\vartheta)\\
		&\preceq_{i_{2}} \mathfrak{p}(\sigma,\varphi)+\mathfrak{p}(\varphi,\vartheta)-\mathfrak{p}(\varphi,\varphi)+s[\mathfrak{d}(\sigma,\varphi)+\mathfrak{d}(\varphi,\vartheta)]\\
		&\preceq_{i_{2}} s[\mathfrak{p}(\sigma,\varphi)+\mathfrak{p}(\varphi,\vartheta)-\mathfrak{p}(\varphi,\varphi)+\mathfrak{d}(\sigma,\varphi)+\mathfrak{d}(\varphi,\vartheta)]\\
		&=s[\delta_{cb}(\sigma,\varphi)+\delta_{cb}(\varphi,\vartheta)-\delta_{cb}(\varphi,\varphi)]\\
		&\preceq_{i_{2}} s[\delta_{cb}(\sigma,\varphi)+\delta_{cb}(\varphi,\vartheta)]-\delta_{cb}(\varphi,\varphi).
		\end{align*} 
		Therefore, $\delta_{cb}(\sigma,\vartheta)\preceq_{i_{2}} s[\delta_{cb}(\sigma,\varphi)+\delta_{cb}(\varphi,\vartheta)]-\delta_{cb}(\varphi,\varphi)$.\\
		So, $(\OE, \delta_{cb})$ is a bicomplex partial b-metric space.  
	\end{proof}
	\begin{prop}
		Let $(\OE,\mathfrak{p})$ be a bicomplex partial metric space, $r\geq 1$, then $(\OE, \delta_{cb})$ is a bicomplex partial b-metric space with coefficient $s=2^{r-1}$,where $\delta_{cb}$ is defined by $\delta_{cb}(\sigma,\vartheta)=[\mathfrak{p}(\sigma,\vartheta)]^r$.
	\end{prop}
	\begin{proof}
		Let $(\OE,\mathfrak{p})$ be a complex partial metric space.\\
		Then $(i)$, $(ii)$, and $(iii)$ of Definition \ref{dt1} are obvious for the function $\delta_{cb}$. \\
		Let $\sigma, \vartheta, \varphi\in \OE$ be arbitrary.\\
		Now,
		\begin{align*}
		\delta_{cb}(\sigma,\vartheta)&=(\mathfrak{p}(\sigma,\vartheta))^r\\
		&=(\mathfrak{p}(\sigma,\varphi)+\mathfrak{p}(\varphi,\vartheta)-(\mathfrak{p}(\varphi,\varphi))^r\\
		&\preceq_{i_{2}}\frac{(\mathfrak{p}(\sigma,\varphi))^r+(\mathfrak{p}(\varphi,\vartheta)-\mathfrak{p}(\varphi,\varphi))^r}{2}2^r\\
		&\preceq_{i_{2}} 2^{r-1}((\mathfrak{p}(\sigma,\vartheta))^r+(\mathfrak{p}(\varphi,\vartheta))^r-\frac{(\mathfrak{p}(\varphi,\varphi))^r}{2^{r-1}})\\
		&=2^{r-1}((\mathfrak{p}(\sigma,\vartheta))^r+(\mathfrak{p}(\varphi,\vartheta))^r)-(\mathfrak{p}(\varphi,\varphi))^r\\
		&=s(\delta_{cb}(\sigma,\varphi)+\delta_{cb}(\varphi,\vartheta))-\delta_{cb}(\varphi,\varphi).
		\end{align*}
		Therefore, $\delta_{cb}(\sigma,\vartheta)\preceq_{i_{2}} s(\delta_{cb}(\sigma,\varphi)+\delta_{cb}(\varphi,\vartheta))-\delta_{cb}(\varphi,\varphi)$.
		So, $(\OE, \delta_{cb})$ is a complex partial b-metric space.
	\end{proof}
	Every bicomplex partial b-metric $\delta_{cb}$ on a non-empty set $\OE$ generates a topology $\tau_{cb}$ on $\OE$ whose base is the family of open $\delta_{cb}$-balls $\mathcal{B}_{\delta_{cb}}(\sigma,\epsilon)$ where $\tau_{cb}=\{\mathcal{B}_{\delta_{cb}}(\sigma,\epsilon):\sigma\in \OE,\epsilon >0\}$ and $\mathcal{B}_{\delta_{cb}}(\sigma,\epsilon)=\{\vartheta\in \OE:\delta_{cb}(\sigma,\vartheta)\prec_{i_{2}}\epsilon+\delta_{cb}(\sigma,\sigma)\}$. Obviously, the topological space $(\OE,\tau_{cb})$ is $\mathcal{T}_{0}$, but need not be $\mathcal{T}_{1}$. 
	\begin{defn}
		Let $(\OE,\delta_{cb})$ be a bicomplex partial b-metric space with coefficient $s$. Let $\{\sigma_{n}\}$ be any sequence in $\OE$ and $\sigma \in \OE$. Then
		\begin{enumerate}
			\item [(i)]  The sequence $\{\sigma_{n}\}$ is said to be convergent with respect to $\wr_{cb}$ and converges to $\sigma$, if $\lim_{n \to \infty}\delta_{cb}(\sigma_{n}, \sigma)=\delta_{cb}(\sigma, \sigma)$.
			\item [(ii)]  The sequence $\{\sigma_{n}\}$ is said to be Cauchy sequence in $(\OE, \delta_{cb})$ if \\$\lim_{n,m \to \infty}\delta_{cb}(\sigma_n, \sigma_m)$ exists and is finite.
			\item [(iii)]  $(\OE,\delta_{cb})$ is said to be a complete bicomplex partial b-metric space if for every Cauchy sequence $\{\sigma_{n}\}$ in $\OE$ there exists $\sigma\in \OE$ such that \\
			$\lim_{n,m \to \infty}\delta_{cb}(\sigma_n, \sigma_m)=\lim_{n \to \infty}\delta_{cb}(\sigma_{n}, \sigma)=\delta_{cb}(\sigma, \sigma)$.
			\item [(iv)]   A mapping $\sqcup:\OE \to \OE$ is said to be continuous at $\sigma_{0}\in \OE$ if for every $\epsilon>0$, there exists $t>0$ such that $\sqcup(B_{\delta_{cb}}(\sigma_{0},t))\subset B_{\delta_{cb}}(\sqcup(\sigma_{0},\epsilon))$.
		\end{enumerate}
	\end{defn}
	\begin{defn}
		Let $(\OE,\delta_{cb})$ be a complex partial b-metric space with coefficient $s\geq 1$. Let $(\OE, \preceq_{i_{2}})$ be a partially ordered set. A pair $(\sqcup, \sqcap)$ of self-maps of $\OE$ is said to be weakly
		increasing if $\sqcup \sigma\preceq_{i_{2}}  \sqcap \sqcup\sigma$ and $\sqcap \sigma\preceq_{i_{2}} \sqcup \sqcap \sigma$ for all $\sigma\in \OE$. If $\sqcup=\sqcap$, then we have $\sqcup \sigma\preceq_{i_{2}} \sqcup^{2}\sigma$
		for all $\sigma\in\OE$ and in this case, we say that $\sqcup$ is weakly increasing mapping.
	\end{defn}
	\begin{defn}
		Let $(\OE,\delta_{cb})$ be a complex partial b-metric space with coefficient $s\geq 1$. A point $\sigma \in \OE$ is said to be common fixed point for the pair of self mappings $(\sqcup, \sqcap)$ on $\OE$ is such that 
		$\sigma=\sqcup \sigma=\sqcap \sigma$.
	\end{defn}
	\begin{example}
		Let $\OE=[0,\infty)$ endowed with bicomplex partial b-metric $\delta_{cb}: X \times X \to \mathbb{C}^{+}_{2}$ with $\delta_{cb}=(\max\{\sigma,\vartheta\})^2+i_{2}(\max\{\sigma,\vartheta\})^2$ $\forall$ $\sigma,\vartheta \in \OE$.
	\end{example}
	
	It is easy to verify that $(\OE, \delta_{cb})$ is a bicomplex partial b-metric space and note that self distance need not be zero, for example $\delta_{cb}(1,1)=1+i_{2}\neq 0$. Now the bicomplex valued b-metric is not induced by $\delta_{cb}$ is follows, $\mathfrak{d}_{\delta_{cb}}(\sigma,\vartheta)=2\delta_{cb}(\sigma,\vartheta)-\delta_{cb}(\sigma,\sigma)-\delta_{cb}(\vartheta,\vartheta)$ without loss of generality suppose $\sigma\geq \vartheta$ then 
	$\mathfrak{d}_{\delta_{cb}}(\sigma,\vartheta)=2[(\max\{\sigma,\vartheta\})^2+i_{2}(\max\{\sigma,\vartheta\})^2]-({\sigma^2+i_{2}\sigma^2})-({\vartheta^2+i_{2}\vartheta^2})$. Therefore, $\mathfrak{d}_{\delta_{cb}}(\sigma,\vartheta)=\sigma^2-\vartheta^2+i_{2}(\sigma^2-\vartheta^2)$.\\
	Therefore, we have the following proposition.
	\begin{prop}
		Every bicomplex partial b-metric $\delta_{cb}$ is not defines bicomplex b-metric $\mathfrak{d}_{\delta_{cb}}$, where $\mathfrak{d}_{\delta_{cb}}(\sigma,\vartheta)=2\delta_{cb}(\sigma,\vartheta)-\delta_{cb}(\sigma,\sigma)-\delta_{cb}(\vartheta,\vartheta)$ $\forall$ $\sigma, \vartheta\in \OE$.
	\end{prop}
	So, we introduce the new notion generalized bicomplex partial b-metric space.
	\begin{defn}
		A generalized bicomplex partial b-metric on a non-empty set $\OE$ is a function $\delta_{cb}:\OE \times \OE\to \mathbb{C}^{+}_{2}$ such that for all $\sigma, \vartheta, \varphi\in \OE$:
		\begin{enumerate}[label=(\roman*)]
			\item $0\preceq \delta_{cb}(\sigma, \sigma) \preceq \delta_{cb}(\sigma,\vartheta)$ (small self-distances)
			\item $\delta_{cb}(\sigma,\vartheta)=\delta_{cb}(\vartheta,\sigma)$ (symmetry)
			\item $\delta_{cb}(\sigma,\sigma)=\delta_{cb}(\sigma,\vartheta)=\delta_{cb}(\vartheta,\vartheta) \Leftrightarrow \sigma=\vartheta$ (equality)
			\item $\exists$ a real number $s\geq 1$ such that $\delta_{cb}(\sigma,\vartheta)\preceq_{i_{2}} s[\delta_{cb}(\sigma,\varphi)+\delta_{cb}(\varphi,\vartheta)-\delta_{cb}(\varphi,\varphi)]+(\frac{1-s}{2})(\delta_{cb}(\sigma,\sigma)+\delta_{cb}(\vartheta,\vartheta))$ (triangularity).
			
		\end{enumerate}
		A generalized bicomplex partial b-metric space is a pair $(\OE,\delta_{cb})$ such that $\OE$ is a non empty set and $\delta_{cb}$ is generalized bicomplex partial b-metric on $\OE$.The number $s$ is called the coefficient of $(\OE,\delta_{cb})$. 
		
	\end{defn}
	\begin{prop}
		Every generalized bicomplex partial b-metric $\delta_{cb}$ is defines bicomplex b-metric $\mathfrak{d}_{\delta_{cb}}$, where $\mathfrak{d}_{\delta_{cb}}(\sigma,\vartheta)=2\delta_{cb}(\sigma,\vartheta)-\delta_{cb}(\sigma,\sigma)-\delta_{cb}(\vartheta,\vartheta)$ $\forall$ $\sigma, \vartheta\in \OE$.
	\end{prop}
	\begin{proof}
		The proof of this proposition is obvious. 
	\end{proof}
	In 2019, Gunaseelan \cite{1aa} proved some fixed point theorems on complex partial b-metric space as follows
	\begin{thm}
		Let $(\OE,\delta_{cb})$ be any complete complex partial b-metric space with coefficient $s\geq 1$ and $\sqcup: \OE\to \OE$ be a
		mapping satisfying:
		\begin{align*}
		\delta_{cb}(\sqcup\sigma, \sqcup\vartheta)\preceq_{i_{2}} \curlywedge [\delta_{cb}(\sigma, \sqcup\sigma)+\delta_{cb}(\vartheta,\sqcup\vartheta)]
		\end{align*}
		for all $\sigma,\vartheta\in \OE$, where $\curlywedge \in [0, \frac{1}{s}]$. Then $\sqcup$ has a unique fixed point $\sigma^{*}\in \OE$ and $\delta_{cb}(\sigma^{*},\sigma^{*})=0$.
	\end{thm}
	We prove some common fixed point theorems on bicomplex partial b-metric space, inspired by his work.
	
	\section{Main Results}
	\begin{thm}\label{t1}
		Let $(\OE, \delta_{cb})$ be a complete bicomplex partial b-metric space with the coefficient $s\geq 1$ and $\sqcap, \sqcup \colon \OE \rightarrow \OE$ be two weakly increasing mappings such that 
		\begin{align*}
		\delta_{cb}(\sqcup \sigma,\sqcap \vartheta) &\preceq_{i_{2}} \dfrac{\curlyvee \delta_{cb}(\sigma,\sqcup \sigma)\delta_{cb}(\vartheta,\sqcap \vartheta)}{\delta_{cb}(\sigma, \vartheta)}+\curlywedge \delta_{cb}(\sigma, \vartheta)
		\end{align*}
		for all $\sigma, \vartheta \in \OE$, $\delta_{cb}(\sigma, \vartheta)\neq 0$ with $\curlyvee\geq 0, \curlyvee+\curlywedge<1$ or $\delta_{cb}(\sqcup \sigma,\sqcap \vartheta)=0$ if $\delta_{cb}(\sigma, \vartheta)=0$  and all elements on the right side are comparable to each other with respect to the partial order $\preceq_{i_{2}}$. If $\sqcup$ or $\sqcap$ is continuous then the pair $(\sqcup, \sqcap)$ have unique common fixed point $\varphi\in\OE$ and $\delta_{cb}(\varphi,\varphi)=0$. 
	\end{thm}
	\begin{proof}
		Let $\sigma_0$ be an arbitrary point in $\OE$ and define a sequence as follows:
		\begin{align*}
		\sigma_{2k+1}&=\sqcup \sigma_{2k}  \\
		\sigma_{2k+2}&=\sqcap \sigma_{2k+1}, k=0,1,2,\ldots.
		\end{align*} 	  
		Since $\sqcup$ and $\sqcap$ are weakly increasing, 
		\begin{align*}
		\sigma_{1}=\sqcup\sigma_{0}&\preceq_{i_{2}} \sqcap\sqcup\sigma_{0}=\sqcap\sigma_{1}=\sigma_{2} \\
		\sigma_{2}=\sqcap\sigma_{1}&\preceq_{i_{2}} \sqcup\sqcap\sigma_{1}=\sqcup\sigma_{2}=\sigma_{3}.
		\end{align*}
		Continuing this way, we have $\sigma_{1}\preceq_{i_{2}} \sigma_{2}\preceq_{i_{2}} \ldots\preceq_{i_{2}}\sigma_{n}\preceq_{i_{2}}\sigma_{n+1}\ldots$. \\
		Assume that $\delta_{cb}(\sigma_{2k},\sigma_{2k+1})>0$ for all $k\in\mathbb{N}$. If not, then $\sigma_{2k}=\sigma_{2k+1}$ for some $k$. Then for all those $k$, $\sigma_{2k}=\sigma_{2k+1}=\sqcup\sigma_{2k}$ and the proof is completed. Assume that $\delta_{cb}(\sigma_{2k},\sigma_{2k+1})>0$ for $k=0,1,2,\ldots.$ As $\sigma_{2k}$ and $\sigma_{2k+1}$ are comparable, so we have 
		\begin{align*}
		\delta_{cb}(\sigma_{2k+1},\sigma_{2k+2})&=\delta_{cb}(\sqcup\sigma_{2k},\sqcap\sigma_{2k+1})\\
		&\preceq_{i_{2}}\curlyvee\dfrac{\delta_{cb}(\sigma_{2k},\sqcup\sigma_{2k})\delta_{cb}(\sigma_{2k+1},\sqcap\sigma_{2k+1})}{\delta_{cb}(\mu_{2k},\mu_{2k+1})}+\curlywedge\delta_{cb}(\sigma_{2k},\sigma_{2k+1}) \\
		&\preceq_{i_{2}} \curlyvee \delta_{cb}(\sigma_{2k+1},\sqcap\sigma_{2k+2})+\curlywedge\delta_{cb}(\sigma_{2k},\sigma_{2k+1}) \\
		\delta_{cb}(\sigma_{2k+1},\sigma_{2k+2})&\preceq_{i_{2}}\dfrac{\curlywedge}{1-\curlyvee}\delta_{cb}(\sigma_{2k},\sigma_{2k+1}).
		\end{align*}
		Now with $\vartheta=\dfrac{\curlywedge}{1-\curlyvee}$, we have
		
		\begin{align*}
		\delta_{cb}(\sigma_{2k+1},\sigma_{2k+2})&\preceq_{i_{2}} h\delta_{cb}(\sigma_{2k},\sigma_{2k+1})\preceq_{i_{2}}\ldots\preceq_{i_{2}} h^{2k+1}\delta_{cb}(\sigma_{0},\sigma_{1}).
		\end{align*}
		For $n>m$, we get 
		\begin{multline*}
		\begin{split}
		\delta_{cb}(\sigma_{m},\sigma_{n})&\preceq_{i_{2}} s\delta_{cb}(\sigma_{m},\sigma_{m+1})+s^{2}\delta_{cb}(\sigma_{m+1},\sigma_{m+2})+\cdots+s^{n}\delta_{cb}(\sigma_{n-1},\sigma_{n})\\
		&-\delta_{cb}(\sigma_{m+1},\sigma_{m+1})-\delta_{cb}(\sigma_{m+2},\sigma_{m+2})-\delta_{cb}(\sigma_{m+3},\sigma_{m+3})\\
		&-\cdots-\delta_{cb}(\sigma_{n-1},\sigma_{n-1})\\
		&\preceq_{i_{2}} (s\vartheta^{m}+s^{2}\vartheta^{m+1}+\cdots+s^{n}\vartheta^{n-1})\delta_{cb}(\sigma_{0},\sigma_{1})\\
		&=s\vartheta^{m}(1+s\vartheta+\cdots +s^{n-1}\vartheta^{n-m-1})\delta_{cb}(\sigma_{0},\sigma_{1})\\
		&\preceq_{i_{2}} \frac{s\vartheta^{m}}{1-s\vartheta}\delta_{cb}(\sigma_{0},\sigma_{1}).
		\end{split}
		\end{multline*}
		Consequently, 
		\begin{align*}
		|\delta_{cb}(\sigma_{m},\sigma_{n})|\leq \dfrac{(s\vartheta)^m}{1-\vartheta}|\delta_{cb}(\sigma_{1},\sigma_{0})|\rightarrow 0
		\end{align*}
		as $m,n\rightarrow \infty$ which implies that $\lim_{m,n\rightarrow \infty}\delta_{cb}(\sigma_{m},\sigma_{n})=0$ such that $\sigma_{n}$ is a Cauchy sequence in $\OE$. Since $(\OE,\delta_{cb})$ is complete, there exists $\varphi\in \OE$ such that $\sigma_{n}\rightarrow \varphi$ and 
		\begin{align*}
		\delta_{cb}(\varphi,\varphi)=\lim_{m,n\rightarrow \infty}\delta_{cb}(\varphi,\sigma_{n})=\lim_{m,n\rightarrow \infty}\delta_{cb}(\sigma_{n},\sigma_{n})=0.
		\end{align*} 
		Without loss of generality, suppose that $\sqcap$ is continuous in  $(\OE,\delta_{cb})$. Therefore $\sqcap\sigma_{2n+1}\rightarrow\sqcap\varphi$ in $(\OE,\delta_{cb})$. That is
		\begin{align*}
		\delta_{cb}(\sqcap\varphi,\sqcap\varphi)=\lim_{n\rightarrow \infty}\delta_{cb}(\sqcap\varphi,\sqcap\sigma_{2n+1})=\lim_{n\rightarrow \infty}\delta_{cb}(\sqcap\sigma_{2n+1},\sqcap\sigma_{2n+1}).
		\end{align*} 
		But 
		\begin{align*}
		\delta_{cb}(\sqcap\varphi,\sqcap\varphi)=\lim_{n\rightarrow \infty}\delta_{cb}(\sqcap\sigma_{2n+1},\sqcap\sigma_{2n+1})=\lim_{n\rightarrow \infty}\delta_{cb}(\sigma_{2n+2},\sigma_{2n+2})=0.
		\end{align*} 
		Next we will prove $\varphi$ is a fixed point of $\sqcap$. 
		\begin{align*}
		\delta_{cb}(\sqcap\varphi,\varphi)\preceq_{i_{2}} s\{\delta_{cb}(\sqcap\varphi,\sqcap\sigma_{2n+1})+\delta_{cb}(\sqcap\sigma_{2n+1},\varphi)\}-\delta_{cb}(\sqcap\sigma_{2n+1},\sqcap\sigma_{2n+1}).
		\end{align*} 
		As $n\rightarrow \infty$, we obtain $|\delta_{cb}(\sqcap\varphi,\varphi)|\leq 0$. Thus, $\delta_{cb}(\sqcap\varphi,\varphi)=0$. Hence $\delta_{cb}(\varphi,\varphi)=\delta_{cb}(\varphi,\sqcap\varphi)=\delta_{cb}(\sqcap\varphi,\sqcap\varphi)=0$ and so $\sqcap\varphi=\varphi$. Therefore $\sqcup\varphi=\sqcap\varphi=\varphi$ and $\delta_{cb}(\varphi,\varphi)=0$. \\
		Finally, we show that $\gamma$ is the unique common fixed point of $\sqcup$ and $\sqcap$. Assume that there
		exists another point $\gamma$ such that $\gamma=\sqcup\gamma=\sqcap\gamma$. Then  $\delta_{cb}(\varphi,\gamma)>0$, we obtain  
		\begin{align*}
		\delta_{cb}(\varphi,\gamma)&=\delta_{cb}(\sqcup\varphi,\sqcap \gamma)\\
		&\preceq_{i_{2}} \curlyvee \dfrac{\delta_{cb}(\varphi,\sqcup\varphi)\delta_{cb}(\gamma,\sqcap \gamma)}{\delta_{cb}(\varphi,\gamma)}+\curlywedge\delta_{cb}(\varphi,\gamma)\\
		&\preceq_{i_{2}} \curlywedge\delta_{cb}(\varphi,\gamma),
		\end{align*}
		which is a contradiction and hence $\varphi=\gamma$. Therefore, $\varphi$ is a
		unique common fixed point of $\sqcup$ and $\sqcap$.
	\end{proof}
	In the absence of the continuity condition for the mapping $\sqcap$, we get the the following Theorem.  
	
	\begin{thm}\label{t2}
		Let $(\OE, \delta_{cb})$ be a complete bicomplex partial b-metric space with the coefficient $s\geq 1$ and $\sqcap, \sqcup \colon \OE \rightarrow \OE$ be two weakly increasing mappings such that 
		\begin{align*}
		\delta_{cb}(\sqcup \sigma,\sqcap \vartheta) &\preceq_{i_{2}} \dfrac{\curlyvee \delta_{cb}(\sigma,\sqcup \sigma)\delta_{cb}(\vartheta,\sqcap \vartheta)}{\delta_{cb}(\sigma, \vartheta)}+\curlywedge \delta_{cb}(\sigma, \vartheta)
		\end{align*}
		for all $\sigma, \vartheta \in \OE$, $\delta_{cb}(\sigma, \vartheta)\neq 0$ with $\curlyvee\geq 0, \curlywedge\geq 0, \curlyvee+\curlywedge<1$ or $\delta_{cb}(\sqcup \sigma,\sqcap \vartheta)=0$ if $\delta_{cb}(\sigma, \vartheta)=0$  and all elements on the right side are comparable to each other with respect to the partial order $\preceq_{i_{2}}$.
		Suppose $\OE$ satisfying the condition that, for every increasing sequence $\{\sigma_{n}\}$ with $\sigma_{n}\rightarrow \varphi$ in $\OE$, we necessarily have $\varphi=\sup \sigma_{n}$, then the pair $(\sqcup, \sqcap)$ have unique common fixed point $\varphi\in \OE$ and $\delta_{cb}(\varphi, \varphi)=0$. 
	\end{thm}
	\begin{proof}
		Given that $\sigma_{n}\preceq_{i_{2}} \varphi$ for all $n\in\mathbb{N}$. Following the proof of the Theorem \ref{t1}, it is enough to prove that $\varphi$ is a fixed point of $\sqcup$. Suppose $\varphi$ is not a fixed point, then we have $\delta_{cb}(\varphi, \sqcup\varphi)=\omega>0$ for some  $\omega\in \mathbb{C}^{+}_{2}$, we obtain 
		\begin{align*}
		\omega&\preceq_{i_{2}} s[\delta_{cb}(\varphi, \sigma_{2n+2})+\delta_{cb}(\sigma_{2n+2},\sqcup\varphi)]-\delta_{cb}(\sigma_{2n+2},\sigma_{2n+2})\\
		&=s[\delta_{cb}(\varphi, \sigma_{2n+2})+\delta_{cb}(\sqcap\sigma_{2n+1},\sqcup\varphi)]-\delta_{cb}(\sigma_{2n+2},\sigma_{2n+2})\\
		&\preceq_{i_{2}} s[\delta_{cb}(\varphi, \sigma_{2n+2})+\curlyvee\dfrac{\delta_{cb}(\sigma_{2n+1},\sqcap\sigma_{2n+1})\delta_{cb}(\varphi,\sqcup\varphi)}{\delta_{cb}(\sigma_{2n+1},\varphi)}\\
		&+\curlywedge\delta_{cb}(\sigma_{2n+1},\varphi)]-\delta_{cb}(\sigma_{2n+2},\sigma_{2n+2}).
		\end{align*}
		Suppose $\delta_{cb}(\varphi, \varphi)=0$, taking limit as $n\rightarrow\infty$, we have $\omega\preceq_{i_{2}} 0$, which is a contradiction. Therefore $\varphi$ is a fixed point of $\sqcup$.\\
		For $\delta_{cb}(\varphi, \varphi)\neq0$, taking limit as $n\rightarrow\infty$, we have $\omega\preceq_{i_{2}} \curlyvee\delta_{cb}(\varphi,\sqcup\varphi)+\curlywedge\delta_{cb}(\varphi, \varphi)$ and so, $|\omega|\leq (\curlyvee+\curlywedge)|\omega|$, since $\curlyvee+\curlywedge<1$, we get a contradiction, which implies that $\varphi=\sqcup\varphi$. Therefore, by Theorem \ref{t1}, we get $\sqcup\varphi=\sqcap\varphi=\varphi$ and $\delta_{cb}(\varphi, \varphi)=0$.\\ 
		Finally, we show that $\gamma$ is the unique common fixed point of $\sqcup$ and $\sqcap$. Assume that there
		exists another point $\gamma$ such that $\gamma=\sqcup\gamma=\sqcap\gamma$. Then  $\delta_{cb}(\varphi,\gamma)>0$, we obtain  
		\begin{align*}
		\delta_{cb}(\varphi,\gamma)&=\delta_{cb}(\sqcup\varphi,\sqcap \gamma)\\
		&\preceq_{i_{2}} \curlyvee \dfrac{\delta_{cb}(\varphi,\sqcup\varphi)\delta_{cb}(\gamma,\sqcap \gamma)}{\delta_{cb}(\varphi,\gamma)}+\curlywedge\delta_{cb}(\varphi,\gamma)\\
		&\preceq_{i_{2}} \curlywedge\delta_{cb}(\varphi,\gamma),
		\end{align*}
		which is a contradiction and hence $\varphi=\gamma$. Therefore, $\varphi$ is a
		unique common fixed point of $\sqcup$ and $\sqcap$.
	\end{proof}
	
	\begin{example}
		Let $\OE=\{1,2,3,4\}$ be endowed with the order $\sigma\preceq_{i_{2}} \vartheta$ if and only if $\vartheta\leq \sigma$. Then $\preceq_{i_{2}}$ is a partial order in $\OE$. Define the bicomplex partial b-metric space $\delta_{cb}:\OE\times \OE\to \mathbb{C}^{+}_{2}$ as follows:
		\begin{center}
			\begin{tabular}{| p{6cm} | l |}
				\hline
				$(\sigma,\vartheta)$ & $\delta_{cb}(\sigma,\vartheta)$\\ 
				\hline
				(1,1), (2,2)& 0 \\ 
				\hline
				(1,2),(2,1),(1,3),(3,1),(2,3),(3,2),(3,3) &$e^{2iy}$ \\ 
				\hline
				(1,4),(4,1),(2,4),(4,2),(3,4),(4,3),(4,4)& $9e^{2iy}$ \\
				\hline
			\end{tabular}.
		\end{center}
		It is easy to verify that $(\OE,\delta_{cb})$ is a complete bicomplex partial b-metric space with the coefficient $s\geq 1$ for $y\in[0,\frac{\pi}{2}]$. Define $\sqcup,\sqcap: \OE\to \OE$ by $\sqcup\sigma=1$,
		\begin{equation*}
		\sqcap(\sigma)=
		\begin{cases}
		1 & \text{if}\ \sigma\in \{1,2,3\} \\
		2 & \text{if}\,\sigma=4
		\end{cases}.
		\end{equation*}
		Then, $\sqcup$ and $\sqcap$ are weakly increasing with respect to $\preceq_{i_{2}}$ and continuous. 
		Now for $\curlyvee=\curlywedge=\frac{1}{9}$, we consider the following cases:
		\begin{enumerate}[label=(\alph*)]
			\item If $\sigma=1$ and $\vartheta\in \OE-\{4\}$, then $\sqcup(\sigma)=\sqcap(\vartheta)=1$ and $\delta_{cb}(\sqcup(\sigma),\sqcap (\vartheta))=0$ and the conditions of Theorem \ref{t1} satisfied.
			\item  If $\sigma=1$, $\vartheta=4$, then $\sqcup\sigma=1$, $\sqcap\vartheta=2$,
			\begin{align*}
			\delta_{cb}(\sqcup\sigma,\sqcap\vartheta)=e^{2iy}&\preceq_{i_{2}} 9\curlywedge e^{i2y}\\
			&=\curlyvee \frac{(0)9e^{2iy}}{9e^{2iy}}+\curlywedge 9e^{2iy}\\
			&=\curlyvee\frac{\delta_{cb}(\sigma,\sqcup\sigma)\delta_{cb}(\vartheta,\sqcap\vartheta)}{\delta_{cb}(\sigma,\vartheta)}+\curlywedge\delta_{cb}(\sigma,\vartheta).
			\end{align*}
			\item  If $\sigma=2$, $\vartheta=4$, then $\sqcup\sigma=1$, $\sqcap\vartheta=2$,
			\begin{align*}
			\delta_{cb}(\sqcup\sigma,\sqcap\vartheta)=e^{2iy}&\preceq_{i_{2}} (\curlyvee +9\curlywedge) e^{i2y}\\
			&=\curlyvee \frac{(e^{2iy})9e^{2iy}}{9e^{2iy}}+\curlywedge 9e^{2iy}\\
			&=\curlyvee\frac{\delta_{cb}(\sigma,\sqcup\sigma)\delta_{cb}(\vartheta,\sqcap\vartheta)}{\delta_{cb}(\sigma,\vartheta)}+\curlywedge\delta_{cb}(\sigma,\vartheta).		
			\end{align*}
			\item  If $\sigma=3$, $\vartheta=4$, then $\sqcup \sigma=1$, $\sqcap \vartheta=2$,
			\begin{align*}
			\delta_{cb}(\sqcup\sigma,\sqcap\vartheta)=e^{2iy}&\preceq_{i_{2}} (\curlyvee +9\curlywedge) e^{i2y}\\
			&=\curlyvee \frac{(e^{2iy})9e^{2iy}}{9e^{2iy}}+\curlywedge 9e^{2iy}\\
			&=\curlyvee\frac{\delta_{cb}(\sigma,\sqcup\sigma)\delta_{cb}(\vartheta,\sqcap\vartheta)}{\delta_{cb}(\sigma,\vartheta)}+\curlywedge\delta_{cb}(\sigma,\vartheta).	
			\end{align*} 
			\item  If $\sigma=4$, $\vartheta=4$, then $\sqcup\sigma=1$, $\sqcap\vartheta=2$,
			\begin{align*}
			\delta_{cb}(\sqcup\sigma,\sqcap\vartheta)=e^{2iy}&\preceq_{i_{2}} 9(\curlyvee +\curlywedge) e^{i2y}\\
			&=\curlyvee \frac{(9e^{2iy})9e^{2iy}}{9e^{2iy}}+\curlywedge 9e^{2iy}\\
			&=\curlyvee\frac{\delta_{cb}(\sigma,\sqcup\sigma)\delta_{cb}(\vartheta,\sqcap\vartheta)}{\delta_{cb}(\sigma,\vartheta)}+\curlywedge\delta_{cb}(\sigma,\vartheta).	
			\end{align*} 
			Moreover for $\curlyvee=\curlywedge=\frac{1}{9}$, with $\curlywedge+\curlywedge=\frac{2}{9}<1$, the conditions of Theorem \ref{t1} are satisfied. Therefore, $1$ is the unique common fixed point of $\sqcup$ and $\sqcap$.
		\end{enumerate}
	\end{example}
	
	\section{Application}
	Consider the following system of Urysohn type integral equations.
	\begin{equation}
	\begin{cases}
	\sigma(q)=b(q)+\int_{x}^{y}G_{1}(q,s,\sigma(s))ds\\
	\sigma(q)=b(q)+\int_{x}^{y}G_{2}(q,s,\sigma(s))ds,\label{cq1}
	\end{cases}       
	\end{equation}
	where 
	\begin{enumerate}[label=(\Alph*)]
		\item $\sigma(q)$ is a unknown variables for each $q\in[x,y]$, $x>0$,
		\item $b(q)$ is the deterministic free term defined for $q\in [x,y]$,
		\item $G_{1}(q,s)$ and $G_{2}(q,s)$ are deterministic kernels defined for $q,s\in [x,y]$.
	\end{enumerate}
	Let $\OE=(C[x,y],\mathbb{R}^{n})$, $q>0$ and $\delta_{cb} :\OE\times \OE\to \mathbb{R}^{n}$ defined by
	\begin{align*}
	\delta_{cb}(\sigma,\vartheta)=|\sigma-\vartheta|^{2}+2+i_{2}(|\sigma-\vartheta|^{2}+2),
	\end{align*}
	for all $\sigma,\vartheta\in \OE$.\\
	Obviously $(C[x,y], \mathbb{R}^{n},\delta_{cb})$ is a complete bicomplex partial b-metric space with the constant $s\geq 1$. Further let us consider a Urysohn type integral system as \eqref{cq1} under the following conditions:
	\begin{enumerate}[start=0,label={(R\arabic*)}]
		\item $b(q)\in \OE$; \label{cq2}
		\item for all $q,s\in [x,y]$, we have\label{cq4}
		\begin{align*}
		G_{1}(q,s,\sigma(s))\preceq_{i_{2}} G_{2}(q,s,b(s)+\int_{x}^{y}G_{1}(q,s,\sigma(s))ds)
		\end{align*}
		and
		\begin{align*}
		G_{2}(q,s,\sigma(s))\preceq_{i_{2}} G_{1}(q,s,b(q)+\int_{x}^{y}G_{2}(q,s,\sigma(s))ds);
		\end{align*}
		\item $G_{1}, G_{2}:[x,y]\times[x,y]\times \mathbb{R}^{n}\to \mathbb{R}^{n}$ are continuous functions satisfying\label{cq3}
		\begin{align*}
		|G_{1}(q,s,\sigma(s))-G_{2}(q,s,\vartheta(s))|\preceq_{i_{2}} \sqrt{\frac{|\sigma-\vartheta|^{2}}{2(y-x)}-\frac{2}{y-x}}.
		\end{align*}
	\end{enumerate}
	\begin{thm}
		Let $(C[x,y],\mathbb{R}^{n},\wp_{cb})$ be a complete bicomplex partial b-metric space, then the system \eqref{cq1} under the conditions \ref{cq2}-\ref{cq3} have a unique common solution.
	\end{thm}  
	\begin{proof}
		For $\sigma,\vartheta\in \OE$ and $q\in [x,y]$, we define the continuous mappings $ 
		\sqcup,\sqcap :\OE\to \OE$ by
		\begin{align*}
		\sqcup \sigma(q)=b(q)+\int_{x}^{y}G_{1}(q,s,\sigma(s))ds, 
		\end{align*}
		and
		\begin{align*}
		\sqcap\sigma(q)=b(q)+\int_{x}^{y}G_{2}(q,s,\sigma(s))ds. 
		\end{align*}
		From the condition \ref{cq4}, the mappings $\sqcup$ and $\sqcap$ are weakly increasing with respect to $\preceq_{i_{2}}$. Indeed, for all $q\in [x,y]$, we have
		\begin{align*}
		\sqcup \sigma(q)&=b(q)+\int_{x}^{y}G_{1}(q,s,\sigma(s))ds\\
		&\preceq_{i_{2}} b(q)+\int_{x}^{y}G_{2}(s,k,b(s)+\int_{x}^{y}G_{1}(q,s,\sigma(s))ds)\\
		&=b(q)+\int_{x}^{y}G_{2}(q,s,\sqcup \sigma(s))ds\\ 
		&=\sqcap(\sqcup\sigma(q) ).
		\end{align*}
		Therefore $\sqcup\sigma(q)\preceq_{i_{2}} \sqcap(\sqcup\sigma(q))$. Similarly, one can easily see that $\sqcap\sigma(q) \preceq_{i_{2}} \sqcup(\sqcap \sigma(q) )$.
		Next we have
		\begin{align*}
		\delta_{cb}(\sqcup \sigma(q),\sqcap \vartheta(q))&=|\sqcup \sigma (q)-\sqcap \vartheta(q)|^{2}+2+i_{2}(|\sqcup \sigma (q)-\sqcap \vartheta(q)|^{2}+2)\\
		&=\int_{x}^{y}|G_{1}(q,s,\sigma(s))-G_{2}(q,s,\vartheta(s))|^{2}ds+2\\
		&+i_{2}\bigg(\int_{x}^{y}|G_{1}(q,s,\sigma(s))-G_{2}(q,s,\vartheta(s))|^{2}ds+2\bigg)\\
		&\preceq_{i_{2}}\int_{x}^{y}\bigg(\frac{|\sigma-\vartheta|^{2}}{2(y-x)}-\frac{2}{y-x}\bigg)ds+2\\
		&+i_{2}\bigg(\int_{x}^{y}\bigg(\frac{|\sigma-\vartheta|^{2}}{2(y-x)}-\frac{2}{y-x}\bigg)ds+2\bigg)\\
		&=\frac{|\sigma-\vartheta|^{2}}{2}+i_{2}\bigg(\frac{|\sigma-\vartheta|^{2}}{2}\bigg)\\
		&\preceq_{i_{2}}\frac{|\sigma-\vartheta|^{2}}{2}+1+i_{2}\bigg(\frac{|\sigma-\vartheta|^{2}}{2}+1\bigg)\\
		&=\curlywedge (|\sigma-\vartheta|^{2}+2+i_{2}(|\sigma-\vartheta|^{2}+2))\\
		&=\curlywedge \delta_{cb}(\sigma,\vartheta).
		\end{align*}
		Hence, all the conditions of Theorem \ref{t1} are satisfied for $\curlyvee+\curlywedge(=\frac{1}{2})<1$ with $\curlyvee=0$. Therefore the system of integral equations \eqref{cq1} has a unique common solution. 
	\end{proof}
	\section{Conclusion}
	In this paper, we proved some common fixed point theorems on bicomplex partial b-metric space under rational type contraction mappings. An illustrative example and application on bicomplex partial b-metric space is given. 
	
\end{document}